\documentclass[12 pt,twoside]{article}
\usepackage{graphicx}

\newtheorem{theorem}{\bf Theorem}[section]

\newtheorem{lemma}[theorem]{\bf Lemma}

{{\footnotesize\sc }}
\input{amssym}

\begin{document}
\title{{\large The Ramsey number of loose paths in 3-uniform hypergraphs}}
\vspace{1cm}\author{\small  L. Maherani$^{\footnotesize\textrm{a}}$, G.R. Omidi$^{\footnotesize\textrm{a},\textrm{b},\textrm{1}}$,
 G. Raeisi$^{\footnotesize\textrm{c}}$, M. Shahsiah$^{\footnotesize\textrm{a}}$\medskip \\
\\
\footnotesize  $^{\textrm{a}}$Department of Mathematical Sciences, Isfahan
University of Technology,\\ \footnotesize Isfahan, 84156-83111, Iran\\
\footnotesize  $^{\textrm{b}}$School of Mathematics, Institute for
Research
in Fundamental Sciences (IPM),\\
\footnotesize  P.O.Box: 19395-5746, Tehran, Iran\\
\footnotesize $^{\textrm{c}}$ Mathematics Department, Faculty of Science, Shahrekord University,
\\ \footnotesize  Shahrekord, P.O. Box 115, Iran\\
\footnotesize{l.maherani@math.iut.ac.ir, romidi@cc.iut.ac.ir}\\
\footnotesize {g.raeisi@sci.sku.ac.ir, m.shahsiah@math.iut.ac.ir}}
\date {}

\footnotesize\maketitle\footnotetext[1] {\tt This research was in part supported
by a grant from IPM (No.90050049)} \vspace*{-.9cm} \footnotesize
\begin{abstract}\rm{}
\footnotesize Recently, asymptotic values of 2-color Ramsey
numbers for loose cycles and also loose paths were determined.
Here we determine the 2-color  Ramsey number of $3$-uniform loose
paths  when one of the paths is significantly larger than the
other:
 for every $n\geq \Big\lfloor\frac{5m}{4}\Big\rfloor$, we show that  $$R(\mathcal{P}^3_n,\mathcal{P}^3_m)=2n+\Big\lfloor\frac{m+1}{2}\Big\rfloor.$$
\\{ {Keywords}:{ \footnotesize Ramsey Number,  Loose Path, Loose Cycle.\medskip}}
\noindent
\\{\footnotesize {AMS Subject Classification}:  05C15, 05C55, 05C65.}
\end{abstract}
%%%%%%%%%%%%%%%%%%%%%%%%%%%%%%%%%%%%%%%%%%%%%%%%%%%%%%%%%%%%%%%%%%%%%%%%%%%%%%%%%%%%%%%%%%%%%%%%%%%%%
%%%%%%%%%%%%%%%%%%%%%%%%%%%%%%%%%%%%%%%%%%%%%%%%%%%%%%%%%%%%%%%%%%%%%%%%%%%%%%%%%%%%%%%%%%%%%%%%%%%%%%
%%%%%%%%%%%%%%%%%%%%%%%%%%%%%%%%%%%%%%%%%%%%%%%%%%%%%%%%%%%%%%%%%%%%%%%%%%%%%%%%%%%%%%%%%%%%%%%%%%%%%%
%%%%%%%%%%%%%%%%%%%%%%%%%%%%%%%%%%%%%%%%%%%%%%%%%%%%%%%%%%%%%%%%%%%%%%%%%%%%%%%%%%%%%%%%%%%%%%%%%%%%%%
%%%%%%%%%%%%%%%%%%%%%%%%%%%%%%%%%%%%%%%%%%%%%%%%%%%%%%%%%%%%%%%%%%%%%%%%%%%%%%%%%%%%%%%%%%%%%%%%%%%%%%
%%%%%%%%%%%%%%%%%%%%%%%%%%%%%%%%%%%%%%%%%%%%%%%%%%%%%%%%%%%%%%%%%%%%%%%%%%%%%%%%%%%%%%%%%%%%%%%%%%%%%%
\small
\section{\normalsize{Introduction}}

A\textit{ hypergraph }$\mathcal{H}$ is a pair $\mathcal{H}=(V,E)$,
where $V$ is a finite
 nonempty set (the set of vertices) and $E$ is a collection of
 distinct nonempty subsets of $V$ (the set of edges). A
 \textit{$k$-uniform hypergraph }  is a hypergraph such that all its
 edges have size $k$. For two $k$-uniform hypergraphs $\mathcal{H}$ and $\mathcal{G}$, the \textit{Ramsey number} $R(\mathcal{H},\mathcal{G})$ is
  the smallest number $N$ such that, in any red-blue
coloring of the edges of the complete $k$-uniform hypergraph
$K^k_N$ on $N$ vertices there is either a red copy of
$\mathcal{H}$ or a blue copy of $\mathcal{G}$. There are several
natural definitions for a cycle and a path in a uniform
hypergraph. Here we consider the one called loose. A $k$-uniform
{\it loose cycle} $\mathcal{C}_n^k$ (shortly, a {\it cycle of
length $n$}),  is a hypergraph with vertex set
$\{v_1,v_2,\ldots,v_{n(k-1)}\}$ and with the set of $n$ edges
$e_i=\{v_1,v_2,\ldots, v_k\}+i(k-1)$, $i=0,1,\ldots, n-1$, where
we use mod $n(k-1)$ arithmetic and adding a number $t$ to a set
$H=\{v_1,v_2,\ldots, v_k\}$ means a shift, i.e. the set obtained
by adding $t$ to subscripts of each element of $H$. Similarly, a
$k$-uniform {\it loose path} $\mathcal{P}_n^k$ (simply, a {\it
path of length $n$}), is a hypergraph with vertex set
$\{v_1,v_2,\ldots,v_{n(k-1)+1}\}$  and with the set of $n$ edges
$e_i=\{v_1,v_2,\ldots, v_k\}+i(k-1)$, $i=0,1,\ldots, n-1$ and we denote this path by $e_0e_1\cdots e_{n-1}$.
 For $k=2$ we get the usual definitions of a cycle and a path. In
 this case, a classical result in graph theory (see \cite{Ramsey number of
paths}) states that
$R(P_n,P_m)=n+\big\lfloor\frac{m+1}{2}\big\rfloor,$ where $n\geq
m\geq 1$. Moreover, the exact values of $R(P_n,C_m)$ and
$R(C_n,C_m)$ for positive integers $n$ and $m$ are
determined\cite{survey}. For $k=3$ it was proved in \cite{Ramsy
number of loose cycle} that $R(\mathcal{C}^3_n, \mathcal{C}^3_n)$,
and consequently $R(\mathcal{P}^3_n, \mathcal{P}^3_n)$ and
$R(\mathcal{P}^3_n, \mathcal{C}^3_n)$, are asymptotically equal to
$\frac{5n}{2}$. Subsequently, Gy\'{a}rf\'{a}s et. al. in
\cite{Ramsy number of loose cycle for k-uniform} extended this
result to the $k$-uniform loose cycles and proved that
$R(\mathcal{C}^k_n, \mathcal{C}^k_n)$, and consequently
$R(\mathcal{P}^k_n, \mathcal{P}^k_n)$ and $R(\mathcal{P}^k_n,
\mathcal{C}^k_n)$, are asymptotically equal  to
$\frac{1}{2}(2k-1)n$. For small cases, Gy\'{a}rf\'{a}s et. al.
(see \cite{subm}) proved that
$R(\mathcal{P}^k_3,\mathcal{P}^k_3)=R(\mathcal{P}^k_3,\mathcal{C}^k_3)=R(\mathcal{C}^k_3,\mathcal{C}^k_3)+1=3k-1$ and
$R(\mathcal{P}_4^k,\mathcal{P}_4^k)=R(\mathcal{P}_4^k,\mathcal{C}_4^k)=R(\mathcal{C}_4^k,\mathcal{C}_4^k)+1=4k-2$. To see a survey on
Ramsey numbers involving cycles see \cite{Rad}.
\bigskip

It is easy to see that $N=(k-1)n+\lfloor\frac{m+1}{2}\rfloor$ is a
lower bound for the Ramsey number
$R(\mathcal{P}^k_n,\mathcal{P}^k_m)$. To show this, partition the
vertex set of $\mathcal{K}_{N-1}^k$ into parts $A$ and $B$, where
$|A|=(k-1)n$ and $|B|=\lfloor\frac{m+1}{2}\rfloor-1$, color all
edges that contain a vertex of $B$ blue, and the rest red. Now,
this coloring can not contain a red copy of $\mathcal{P}_n^k$,
since such a copy has $(k-1)n+1$ vertices. Clearly the longest
blue path has length at most $m-1$, which proves our claim. Using
the same argument we can see that $N$ and $N-1$ are the lower
bounds for $R(\mathcal{P}^k_n,\mathcal{C}^k_m)$ and
$R(\mathcal{C}^k_n,\mathcal{C}^k_m)$, respectively. In
\cite{subm}, motivated by the above facts and some other results,
the authors conjectured that these lower bounds give the exact
values of the mentioned Ramsey numbers for $k=3$. In this paper,
we consider this problem and we prove that
$R(\mathcal{P}^3_n,\mathcal{P}^3_m)=2n+\Big\lfloor\frac{m+1}{2}\Big\rfloor$
for every $n\geq \lfloor\frac{5m}{4}\rfloor$. Throughout the
paper, for a 2-edge coloring of a uniform hypergraph
$\mathcal{H}$, say red and blue, we denote by $\mathcal{F}_{red}$
and $\mathcal{F}_{blue}$ the induced hypergraph on edges of colors
red and blue,
respectively.\\

%%%%%%%%%%%%%%%%%%%%%%%%%%%%%%%%%%%%%%%%%%%%%%%%%%%%%%%%%%%%%%%%%%%%%%%%%%%%%%%%%%%%%%%%%%%%%%%%%%%%%%
%%%%%%%%%%%%%%%%%%%%%%%%%%%%%%%%%%%%%%%%%%%%%%%%%%%%%%%%%%%%%%%%%%%%%%%%%%%%%%%%%%%%%%%%%%%%%%%%%%%%%%
\medskip
\section{\normalsize Preliminaries}

\medskip
In this section, we present some lemmas which are essential in the
proof of the main results.

\medskip
\begin{lemma}\label{No Cn}
Let $n\geq m\geq 3$ and
$\mathcal{K}^k_{(k-1)n+\lfloor\frac{m+1}{2}\rfloor}$ be 2-edge
colored red and blue. If $\mathcal{C}_{n}^k\subseteq
\mathcal{F}_{red}$, then either $\mathcal{P}_n^k\subseteq
\mathcal{F}_{red}$ or $\mathcal{P}_m^k\subseteq
\mathcal{F}_{blue}$.
\end{lemma}
%%%%%%%%%%%%%%%%%%%%%%%%%%%%
%%%%%%%%%%%%%%%%%%%%%%%%%%%%
\noindent\textbf{Proof. }Let $e_i=\{v_1,v_2,\ldots, v_k\}+i(k-1)$
(mod $n(k-1)$), $i=0,1,\ldots, n-1$, be the edges of
$\mathcal{C}_{n}^k\subseteq \mathcal{F}_{red}$ and
$W=\{x_1,x_2,\ldots,x_{\lfloor\frac{m+1}{2}\rfloor}\}$
 be the set of the remaining vertices. Set $e_{0}^{\prime}=(e_{0}\setminus\{v_{1}\})\cup\{x_{1}\}$ and for $1\leq i\leq m-1$ let

$$
%\begin{eqnarray*}
e_i^{\prime}= \left\lbrace
\begin{array}{ll}
(e_i\setminus \{v_{i(k-1)+1}\})\cup\{x_{\frac{i+1}{2}}\}  &\mbox{if~}  i~\mbox{is~odd},\vspace{.5 cm}\\
(e_i\setminus\{v_{(i+1)(k-1)+1}\})\cup\{x_{\frac{i+2}{2}}\}
&\mbox{if~}  i~\mbox{is~even}.
\end{array}
\right.\vspace{.2 cm}
$$

\noindent If one of $e_i^{\prime}$ is red,  we have a
monochromatic $\mathcal{P}_n^k\subseteq \mathcal{F}_{red}$,
otherwise $e_0^{\prime}e_1^{\prime}\ldots e^{\prime}_{m-1}$ form a
blue $\mathcal{P}_m^k$, which completes the proof.

$\hfill \blacksquare$
%%%%%%%%%%%%%%%%%%%%%%%%%%%%%%%%%%%%%%%%%%%%%%%%%%%%%%%%%%%%%%%%%%%%%%%%%%%%%%%%%%%%%%%%%%%%%%%%%%%%%%
%%%%%%%%%%%%%%%%%%%%%%%%%%%%%%%%%%%%%%%%%%%%%%%%%%%%%%%%%%%%%%%%%%%%%%%%%%%%%%%%%%%%%%%%%%%%%%%%%%%%%%

\bigskip
Let $\mathcal{P}$ be a loose path and $x,y$ be vertices which are
not in $\mathcal{P}$. By a {\it
$\varpi_{\{v_i,v_j,v_k\}}$-configuration}, we mean a copy of
$\mathcal{P}^3_2$ with edges $\{x,v_i,v_j\}$ and \{$v_j,v_k,y\}$
so that $v_l$'s, $l\in\{i,j,k\}$, belong to two consecutive
edges of $\mathcal{P}$. The vertices $x$ and $y$ are called
the end vertices of this configuration. Using this notation, we
have the following lemmas.

\bigskip
\begin{lemma}\label{spacial configuration}
Let $n\geq 10$,  $\mathcal{K}^3_{n}$ be 2-edge colored red and
blue and $\mathcal{P}$, say in $\mathcal{F}_{red}$, be a maximum
path. Let $A$ be the set of five consecutive vertices of
$\mathcal{P}$.
% with two consecutive edges $e,e'$.
 If
$W=\{x_1,x_2,x_3\}$ is disjoint from $\mathcal{P}$, then we have a
$\varpi_S$-configuration in $\mathcal{F}_{blue}$ with two end
vertices in $W$ and $S\subseteq A$.
\end{lemma}
%%%%%%%%%%%%%%%%%%%%%%%%%%%%
%%%%%%%%%%%%%%%%%%%%%%%%%%%%
\medskip
\noindent\textbf{Proof. }First let $A=e\cup e'$ for two edges
$e=\{v_1,v_2,v_3\}$ and $e'=\{v_3,v_4,v_5\}$. Since
$\mathcal{P}\subseteq \mathcal{F}_{red}$ is maximal, at least one
of the edges $e_1=\{x_1,v_1,v_2\}$ and $e_2=\{v_2,v_3,x_2\}$ must
be blue. If both are blue, then $e_1e_2$ is such a configuration.
So first let $e_1$ be blue and $e_2$ be red. Maximality of
$\mathcal{P}$ implies that at least one of the edges
$e_3=\{x_2,v_1,v_4\}$ or $e_4=\{x_3,v_2,v_5\}$ is blue (otherwise,
replacing $ee'$ by $e_3e_2e_4$ in  $\mathcal{P}$ yields a red path
greater than $\mathcal{P}$, a contradiction), and clearly in each
case we have a $\varpi_S$-configuration. Now, let $e_1$ be red and
$e_2$ be blue. Clearly $e_5=\{v_2,v_4,x_3\}$ is blue and $e_2e_5$
form a $\varpi_S$-configuration. Now let
$A=\{v_1,v_2,\ldots,v_5\}$ where $e_1=\{x,v_1,v_2\}$,
$e_2=\{v_2,v_3,v_4\}$ and $e_3=\{v_4,v_5,y\}$ are three
consecutive edges  of $\mathcal{P}$. If $\{x_i,v_2,v_3\}$ is a red
edge for some $i\in \{1,2,3\}$, then $\{v_3,v_4,x_j\}$ and
$\{v_3,v_5,x_j\}$ are blue for $j\neq i$ and so we are done. By
the same argument the theorem is true if $\{x_i,v_3,v_4\}$ is red.
Now we may assume $\{v_2,v_3,x_i\}$ and $\{v_3,v_4,x_i\}$ are blue
for each $i\in \{1,2,3\}$ and so there is nothing to prove.

$\hfill \blacksquare$
%%%%%%%%%%%%%%%%%%%%%%%%%%%%%%%%%%%%%%%%%%%%%%%%%%%%%%%%%%%%%%%%%%%%%%%%%%%%%%%%%%%%%%%%%%%%%%%%%%%%%%%%%%%
%%%%%%%%%%%%%%%%%%%%%%%%%%%%%%%%%%%%%%%%%%%%%%%%%%%%%%%%%%%%%%%%%%%%%%%%%%%%%%%%%%%%%%%%%%%%%%%%%%%%%%%%%
%%%%%%%%%%%%%%%%%%%%%%%%%%%%%%%%%%%%%%%%%%%%%%%%%%%%%%%%%%%%%%%%%%%%%%%%%%%%%%%%%%%%%%%%%%%%%%%%%%%%%%
%$$$$$$$$$$$$$$$$$$$$$$$$$$$$$$$$$$$$$$$$$$$$$$$$$$$$$$$$$$$$$$$$$$$$$$$$$$$$
\bigskip
\begin{lemma}\label{P5 or P4}
Assume that $n\geq \Big\lfloor\frac{5m}{4}\Big\rfloor$ and
$\mathcal{K}^3_{2n+\lfloor\frac{m+1}{2}\rfloor}$ is 2-edge colored
red and blue. If $\mathcal{P}\subseteq \mathcal{F}_{blue}$ is a
maximum path and  $W$ is the set of the vertices
 which are not covered by $\mathcal{P}$,
then for every 4 consecutive edges $e_1,e_2,e_3,e_4$  of
$\mathcal{P}$  either there is a  $\mathcal{P}^3_{5}\subseteq
\mathcal{F}_{red}$, say $Q$, between $\lbrace e_1,e_2,e_3,e_4
\rbrace$ and $W$ with end vertices in $W$ and with no the last
vertex of $e_4$ as a vertex  such that $\vert W\cap V(Q)\vert \leq
5$ or there is  a $\mathcal{P}^3_{4}\subseteq \mathcal{F}_{red}$,
say $Q$, between $\lbrace e_1,e_2,e_3 \rbrace$ and $W$ with end
vertices in $W$ and with no the last vertex of $e_3$ as a vertex
such
 that $\vert W\cap V(Q) \vert \leq 4$. In  each of the above cases, each vertex of $W$ except one vertex can be considered as the end vertex of $Q$.
\end{lemma}

\noindent\textbf{Proof. }Suppose that $e_1$,$e_2$,$e_3$,$e_4$ be
four consecutive edges in $\mathcal{P}$. Let
 $e_i=\{v_{2i-1},v_{2i},v_{2i+1} \}$, $1\leq i\leq 4$, and $W=\{x_1,...,x_t\}$ and  $T=\{ 1,2,\cdots ,t\}$.
\bigskip

\noindent \textbf{Case 1. }For every $1\leq i,j\leq t$,
$\{v_1,v_2,x_i \}$ and $\{v_2,v_3,x_j \}$ are red.

\medskip
{\it Subcase 1}. For every $1\leq  k,l\leq  t$, the edges
$\{v_{3},v_{4},x_{k}\}$ and $\{v_{4},v_{5},x_{l}\}$ are red.

\medskip  So for each $\{ i_1, i_2, i_3, i_4 \} \in P_{4}(T)$, the edges,
 $\{x_{i_{1}},v_1,v_2\}$,$\{v_2,x_{i_{2}},v_3\}$,$\{v_3,x_{i_{3}},v_4\}$,$\{v_4,v_5,x_{i_{4}}\} $ make
  a red $\mathcal{P}^3_{4}$  with end vertices $x_{i_{1}}$ and $x_{i_{4}}$.

\medskip
{\it Subcase 2}. There exists $1\leq  k\leq  t$, such that the
edge  $\{v_{3},v_{4},x_{k}\}$ is blue.

\medskip So for each $\{ i_1, i_2, i_3, \} \in P_{3}(T)$ with $k \neq i_2 ,i_3$, $\{x_{i_{1}}
,v_1,v_2 \}$,$\{v_2,v_3,x_{i_{2}}
 \}$,$\{x_{i_{2}},v_5,v_4 \}$,\\$\{v_4,v_6,x_{i_{3}}\}$ are the edges of a red desired $\mathcal{P}^3_{4}$ with end vertices $x_{i_{1}}$ and $x_{i_{3}}$.

\medskip
 {\it Subcase 3}. There exists $1\leq  k\leq  t$, such
that  the edge  $\{v_{4},v_{5},x_{k}\}$ is blue.

\medskip  If  for every $1\leq  i,j\leq  t$, the edges  $\{v_{5},v_{6},x_{i}\}$ and $\{v_{6},v_{7},x_{j}\}$ are red,
 then for every $\{ i_1, i_2, i_3, i_4 \}\in P_{4}(T)$ with $i_3 \neq k$,  we can
find a red copy of $\mathcal{P}^3_{5}$ with edges $\{x_{i_{1}},v_1,v_2 \}$,$\{v_2,x_{i_{2}},v_3 \}$,\\$\{v_3,v_4,x_{i_{3}} \}$,$\{x_{i_{3}},v_5,v_6 \}$,$\{v_6,v_7,x_{i_{4}}\}$ and
 end vertices $x_{i_{1}}$ and $x_{i_{4}}$. Otherwise there exists $1\leq   l\leq  t$, such that either
 $\{v_{5,}v_{6},x_{l}\}$ or $\{v_{6},v_{7},x_{l}\}$ is blue. For the first one,  for every $\{ i_1, i_2, i_3, i_4 \}\in P_{4}(T)$ with
  $i_3 \neq k,l$ and $i_4 \neq l$, $\{x_{i_{1}},v_{1},v_{2}\}$,$\{v_{2},x_{i_{2}},v_{3}\}$,$\{v_{3},v_{4},x_{i_{3}}\}$,$\{x_{i_{3}},v_{7},v_6 \}$,$\{v_6,v_8,x_{i_{4}} \}$
   make a red copy of $\mathcal{P}^3_{5}$ with end vertices $x_{i_{1}}$ and $x_{i_{4}}$ and for the second
 one, for every $\{ i_1, i_2, i_3, \}\in P_{3}(T)$ with $l\neq i_2, i_3$  the edges,
    $\{x_{i_{1}},v_{1},v_{2},\},\{v_{2},v_{3},x_{i_{2}}\},\{x_{i_{2}},v_{6},v_{5}\},\{v_{5},x_{i_{3}},x_l \}$ make
     a red $\mathcal{P}^3_{4}$ with end vertices $x_{i_{1}}$ and $y$  where  $y\in \{ x_{i_{3}},x_l \}$.

\bigskip

\noindent \textbf{Case 2. } For some $1\leq i\leq t$,
$\{v_{1},v_{2},x_{i}\}$ is blue.

\medskip
  {\it Subcase 1}. For every $1\leq k,l\leq t$, the edges $\{v_{5},v_{6},x_{k}\}$ and $\{v_{6},v_{7},x_{l}\}$ are red.

\medskip  For each  $\{ i_1, i_2, i_3, i_4 \}\in P_{4}(T)$ with $i_j \neq i$, $1\leq j\leq 4$, the edges,
  $\{ x_{i_{1}},x_{i},v_{3}\}$,$\{v_{3}, x_{i_{2}},v_{2}\}$,\\$\{v_{2},v_{4}, x_{i_{3}}\}$,$\{ x_{i_{3}},v_{5},v_6\}$,$\{v_6,v_7, x_{i_{4}} \}$
   make a red $\mathcal{P}^3_{5}$ with end vertices $y$, $y\in \{ x_{i_{1}},x_i \}$, and $x_{i_{4}}$.

\medskip
  {\it Subcase 2}.  For some $1\leq k\leq t$, $\{v_5,v_6,x_k \}$ is blue.

\medskip  In this case,  for each  $\{ i_1, i_2, i_3, i_4 \}\in P_{4}(T)$ with $i_j \neq i$,   $1\leq j\leq 4$,
 and $i_3, i_4 \neq k$,  the edges  $\{x_{i_{1}},x_{i},v_{3}\}$,$\{v_{3},x_{i_{2}},v_{2}\}$,$\{v_{2},v_{4},x_{i_{3}}\}$,$\{x_{i_{3}},v_{7},v_6 \}$,
 $\{v_6,v_8,x_{i_{4}}\}$ make a red  $\mathcal{P}^3_{5}$ with end vertices $y$, $y\in \{ x_{i_{1}},x_i \}$, and $x_{i_{4}}$.

\medskip
 {\it Subcase 3}.  For some $1\leq k\leq t$, $\{v_6,v_7,x_k \}$ is blue.

 \medskip  In this case, for each  $\{ i_1, i_2, i_3 \}\in P_{3}(T)$ with $i_j \neq i$,  $1\leq j\leq 3$, and $i_2, i_3 \neq k$,  the edges
   $\{x_{i_{1}},x_{i},v_{3}\}$,$\{v_3,v_{2},x_{i_{2}}\}$,$\{x_{i_{2}},v_{4},v_{6}\}$,$\{v_6,v_5,x_{i_{3}}\}$ make a red
    $\mathcal{P}^3_{4}$ with end vertices $y$, $y\in \{ x_{i_{1}},x_i \}$, and $x_{i_3}$.
%%%%%%#######################################################################################3

\bigskip

\noindent \textbf{Case 3. } For some $1\leq i\leq t$,
$\{v_{2},v_{3,}x_{i}\}$ is blue.

\medskip
  {\it Subcase 1}. For every $1\leq k,l\leq t$, the edges $\{v_{3},v_{4},x_{k}\}$ and $\{v_{4},v_{5},x_{l}\}$ are red.

\medskip  For each  $\{ i_1, i_2, i_3 \}\in P_{3}(T)$ with $i_j \neq i$,  $1\leq j\leq 3$,
 $\{x_{i_{1}},x_i,v_1 \}$,$\{v_1,v_2,x_{i_{2}} \}$,$\{x_{i_{2}},v_3,v_4 \}$,\\$\{v_4,v_5,x_{i_{3}} \}$ are
 the edges of a red $\mathcal{P}^3_{4}$ with end vertices  $y$, $y\in \{ x_{i_{1}},x_i \}$, and $x_{i_{3}}$.

\medskip
  {\it Subcase 2}.  For some $1\leq k\leq t$, $\{v_3,v_4,x_k \}$ is blue.

\medskip  In this case,  for each  $\{ i_1, i_2, i_3 \}\in P_{3}(T)$ with $i_j \neq i$,  $1\leq j\leq 3$, and $i_2, i_3 \neq k$,  the edges,
 $\{x_{i_{1}},x_i,v_1 \}$,$\{v_1,v_2,x_{i_{2}} \}$,$\{x_{i_{2}},v_5,v_4 \}$,$\{v_4,v_6,x_{i_{3}} \}$ make
  a red copy of  $\mathcal{P}^3_{4}$ with end vertices $y$, $y\in \{ x_{i_{1}},x_i \}$, and $x_{i_3}$.

\medskip
 {\it Subcase 3}.  For some $1\leq k\leq t$, $\{v_4,v_5,x_k \}$ is blue.

 \medskip If  for every $1\leq  l,h\leq  t$, the edges  $\{v_{5},v_{6},x_{l}\}$ and $\{v_{6},v_{7},x_{h}\}$ are red, then
  for each  $\{ i_1, i_2, i_3, i_4 \}\in P_{4}(T)$ with $i_j \neq i$,  $1\leq j\leq 4$, and $i_3\neq k$,  the edges,
 $\{x_{i_{1}},x_i,v_1 \}$,$\{v_1,x_{i_{2}},v_2 \}$,\\$\{v_2,v_4,x_{i_{3}} \}$,$\{x_{i_{3}},v_5,v_6
\},\{v_6,v_7,x_{i_{4}} \}$
 make a red $\mathcal{P}^3_{5}$ with end vertices $y$, $y\in \{ x_{i_{1}},x_i \}$, and $x_{i_{4}}$.
 Otherwise there exists $1\leq  l\leq  t$, such that
either
 $\{v_{5,}v_{6},x_{l}\}$ or $\{v_{6},v_{7},x_{l}\}$ is blue. For the first one, for each  $\{ i_1, i_2, i_3, i_4 \}\in P_{4}(T)$ with $i_j \neq i$,
  $1\leq j\leq 4$, $i_3\neq k,l$ and $i_4\neq l$,  the edges  $\{x_{i_{1}},x_{i},v_{1}\}$,$\{v_{1},x_{i_{2}},v_{2}\}$,$\{v_{2},v_{4},x_{i_{3}}\}$,
 $\{x_{i_{3}},v_{7},v_6 \}$,$\{v_6,v_8,x_{i_{4}}\}$ make a red copy of $\mathcal{P}^3_{5}$  with end vertices  $y$, $y\in \{ x_{i_{1}},x_i \}$ and
  $x_{i_4}$.  For the second
 one, for every  $\{ i_1, i_2, i_3 \}\in P_{3}(T)$ with $i_j \neq i$,  $1\leq j\leq 3$, and $i_2, i_3 \neq l$,
 $\{\{x_{i_{1}},x_{i},v_{1}\},\{v_{1},v_{2},x_{i_{2}}\},\{x_{i_{2}},v_{4},x_{6}\},
 \{v_{6},v_{5},x_{i_{3}} \}\}$ is the set of the edges of a red $\mathcal{P}^3_{4}$  with end vertices
  $y$, $y\in \{ x_{i_{1}},x_i \}$, and $x_{i_3}$. These observations complete the proof.
  $\hfill \blacksquare$

%$$$$$$$$$$$$$$$$$$$$$$$$$$$$$$$$$$$$$$$$$$$$$$$$$$$$$$$$$$$$$$4
%\medskip%%%%%%%%%%%5

%%%%%%%%%%%%%%%%%%%%%%%%%%%%%%%%%%%%%%%%%%%%%%%%%%%%%%%%%%%%%%%%%%%%%%%%%%%%%%%%%%%%%%%%%%%%%%%%%%%%%%
%%%%%%%%%%%%%%%%%%%%%%%%%%%%%%%%%%%%%%%%%%%%%%%%%%%%%%%%%%%%%%%%%%%%%%%%%%%%%%%%%%%%%%%%%%%%%%%%%%%%%%
\section{\normalsize Main Results}
In this section, we prove that
$R(\mathcal{P}^3_n,\mathcal{P}^3_m)=2n+\Big\lfloor\frac{m+1}{2}\Big\rfloor$
for every $n\geq\lfloor\frac{5m}{4}\rfloor$. First we present
several lemmas which will be our main tools in establishing the main theorem.
\bigskip
\begin{lemma}\label{pm-1 implies pn-1}
Assume that $n= \Big\lfloor\frac{5m}{4}\Big\rfloor$ and
$\mathcal{K}^3_{2n+\lfloor\frac{m+1}{2}\rfloor}$ is 2-edge colored
red and blue. If $\mathcal{P}=\mathcal{P}^3_{m-1}$ is a maximum
blue path, then $\mathcal{P}^3_{n-1}\subseteq \mathcal{F}_{red}$.
\end{lemma}
%%%%%%%%%%%%%%%%%%%%%%%%%%%%
%%%%%%%%%%%%%%%%%%%%%%%%%%%%
\noindent\textbf{Proof. }Let $t=2n+\lfloor\frac{m+1}{2}\rfloor$
and $\mathcal{P}=e_1e_2\ldots e_{m-1}$ be a copy of
$\mathcal{P}_{m-1}^3\subseteq \mathcal{F}_{blue}$ with  edges
$e_i=\{v_1,v_2,v_3\}+2(i-1)$, $i=1,\ldots, m-1$. Set
 $W=V(\mathcal{K}^3_t)\setminus V(\mathcal{P})$. Using Lemma
 \ref{P5 or P4} there is a red path $Q_1$ with end vertices $x_1$ and $y_1$ in $W_1=W$ between $E'_1$ and $W_1$
 where $E_1=\{e_i: i_1=1\leq i \leq 4\}$, $\bar{E}_1=E_1\setminus
 \{e_4\}$ and $E'_1\in \{E_1,\bar{E}_1\}$. Set $i_2=\min\{j: j\in\{i_1+3,i_1+4\},
 e_j\not\in E'_1\}$, $E_2=\{e_i: i_2\leq i\leq i_2+3\}$ and
$\bar{E}_2=E_2\setminus
 \{e_{i_2+3}\}$ and $W_2=(W\setminus V(Q))\cup \{x_1,y_1\}$. Again
 using Lemma \ref{P5 or P4} there is a red path $Q_2$ between
 $E'_2$ and $W_2$ such that $Q_1\cup Q_2$ is a red path with end
 vertices $x_2, y_2$ in $W_2$ where $E'_2\in \{E_2, \bar{E}_2\}$
 and again set $i_3=\min\{j: j\in\{i_2+3,i_2+4\}, e_j\not\in E'_2\}$, $E_3=\{e_i: i_3\leq i\leq
 i_3+3\}$, $\bar{E}_3=E_3\setminus
 \{e_{i_3+3}\}$ and $W_3=(W\setminus V(Q_1\cup Q_2))\cup
 \{x_2,y_2\}$. Since $|W|\geq m$, using Lemma \ref{P5 or P4} by
 continuing the above process we can partition
 $E(\mathcal{P})\setminus\{e_{m-1}\}$ into classes $E'_i$th,
 $|E'_i|\in\{3,4\}$ and at most one class of size  $r\leq 3$ of the last edges such that for
 each $i$, there is a red $Q_i=\mathcal{P}^3_5$ (resp.
 $Q_i=\mathcal{P}^3_4$) between $E'_i$ and $W$ with the
 properties in Lemma \ref{P5 or P4} if $|E'_i|=4$ (resp.
 $|E'_i|=3$) and $\mathcal{P}'=\Large\cup Q_i$ is a red path with end
 vertices $x, y$ in $W$. Let $l_1=|\{i: \mid E'_i\mid=4\}|$ and
$l_2=|\{i: \mid E'_i\mid=3\}|$. So $m-2=4l_1+3l_2+r$, $0\leq r\leq
3$ and  $\mathcal{P}^{\prime}$ has $5l_1+4l_2$ edges. One can
easily check that $5l_1+4l_2\geq \frac{5}{4}(m-2-r)$. Also we have
          $$|W\cap V(\mathcal{P}^{\prime})|\leq 4l_1+3l_2+1=m-1-r.$$
Let $T=V(\mathcal{K}^3_t)\setminus (V(\mathcal{P})\cup
V(\mathcal{P}'))$ and suppose that $m=4k+p$ for some $p$, $0\leq
p\leq 4$. Therefore
 $|T|\geq r+2$ if $p=0,1$ and $|T|\geq r+1$ if $p=2,3$. Now we consider the following cases.

\bigskip
\noindent \textbf{Case 1. }$r=0$.

\medskip
Clearly $|T|\geq 1$ and it is easy to see that
$\mathcal{P}^{\prime}$ contains at least $n-2$ edges. Let
$\{u\}\subseteq T$. The maximality of $\mathcal{P}$ implies that
the edge $e=\{v_{2m-1},x,u\}$ is red and hence
$\mathcal{P}^{\prime}\cup \{e\}$ is  a red copy of
$\mathcal{P}^3_{n-1}$.

\bigskip
\noindent \textbf{Case 2. }$r=1$.

\medskip
In this case, $|T|\geq 2$  and it is easy to see that
$\mathcal{P}^{\prime}$ contains at least $n-3$ edges. Let
$\{u,v\}\subseteq T$. Clearly $\mathcal{P}^{\prime}\cup
\{\{v_{2m-2},x,u\}, \{v_{2m-1},u,v\}\}$ is  a red copy of
$\mathcal{P}^3_{n-1}$.

\bigskip
\noindent \textbf{Case 3. }$r=2$.

\medskip
It is easy to see that  $|T|\geq 3$  and $\mathcal{P}^{\prime}$
contains at least $n-5$ edges. Let $T'=\{u,v,w\}\subseteq T$.
Since $V(\mathcal{P}^{\prime})\cap V(e_{m-3}\cup
e_{m-2})=\emptyset$ by lemma \ref{spacial configuration} there is
a red $\varpi_S$-configuration with $S\subset e_{m-3}\cup e_{m-2}$
and its end vertices in $T'$, say $u$ and $v$. The maximality of
$\mathcal{P}$ implies that the edges $\{v_{2m-2},x,u\}$ and
$\{v_{2m-1},v,w\}$ are red and clearly we
have a red $\mathcal{P}^3_{n-1}$.\\

\bigskip
\noindent \textbf{Case 4. }$r=3$.

\medskip
In this case, for $p\in \{2,3\}$ we have  $|T|\geq 4$  and
$\mathcal{P}^{\prime}$ contains at least $n-5$ edges. Using an
argument similar to case 3 we can complete the proof. Now let
$p\in \{0,1\}$. Then $|T|\geq 5$  and $\mathcal{P}^{\prime}$
contains at least $n-6$ edges. Set $T'=\{u,v,w,z,t\}\subseteq T$.
By Lemma \ref{spacial configuration}, there is a
$\varpi_S$-configuration $C$ with $S\subseteq V(e_{m-3}\cup
e_{m-2})$ and end vertices in $T'$, say $u$ and $v$. Clearly
$\mathcal{P}^{\prime}\cup\{\{y,w,v_{2m-2}\},\{v_{2m-2},z,t\},\{v_{2m-1},t,u\}\}\cup
C$ is a red $\mathcal{P}^3_{n-1}$. These observations complete the
proof. $\hfill \blacksquare$
%%%%%%%%%%%%%%%%%%%%%%%%%%%%%%%%%%%%%%%%%%%%%%%%%%%%%%%%%%%%%%%%%%%%%%%%%%%%%%%%%%%%%%%%%%%%%%%%%%%%%%
%%%%%%%%%%%%%%%%%%%%%%%%%%%%%%%%%%%%%%%%%%%%%%%%%%%%%%%%%%%%%%%%%%%%%%%%%%%%%%%%%%%%%%%%%%%%%%%%%%%%%%

\bigskip
\begin{lemma}\label{a}
Let $n\geq \Big\lfloor\frac{5m}{4}\Big\rfloor$ and
$\mathcal{K}^3_{2n+\lfloor\frac{m+1}{2}\rfloor}$ be 2-edge colored
red and blue. If $\mathcal{P}^3_{n-1}\subseteq \mathcal{F}_{red}$
be a maximum path, then $\mathcal{P}^3_{m}\subseteq
\mathcal{F}_{blue}$.
\end{lemma}
%%%%%%%%%%%%%%%%%%%%%%%%%%%%
%%%%%%%%%%%%%%%%%%%%%%%%%%%%
\noindent\textbf{Proof. }Let $t=2n+\lfloor\frac{m+1}{2}\rfloor$
and $\mathcal{P}=e_1e_2\ldots e_{n-1}$ be a copy of
$\mathcal{P}_{n-1}^3\subseteq \mathcal{F}_{red}$ with end edges
$e_1=\{v_1,v_2,v_3\}$ and
$e_{n-1}=\{v_{2n-3},v_{2n-2},v_{2n-1}\}$. By Lemma \ref{No Cn}, we
 may assume that the subhypergraph induced by $V(\mathcal{P})$ does not have a red copy of $\mathcal{C}^3_n$. Let $W=V(\mathcal{K}^3_t)\setminus V(\mathcal{P})$ and let $2n-2=5q+h$ where $0\leq h<5$.
  Partition the set $V(\mathcal{P})\setminus \{v_1\}$ into $q$ classes $A_1,A_2,\ldots,A_q$ of size five and one class
$A_{q+1}=\{v_{2n-h},\ldots,v_{2n-2},v_{2n-1}\}$ of size $h$ if $h>0$,  so that each class contains consecutive vertices
   of $\mathcal{P}$. Using Lemma \ref{spacial configuration}, there is a blue
$\varpi_{S_1}$-configuration, $\bar{c}_1$, with the set of end
vertices $E_1\subseteq W$ and $S_1\subseteq A_1$. Let $x_1\in E_1$
and $B_1$ be a 2-subset of $W\setminus E_1$. Again by Lemma
\ref{spacial configuration}, there is a blue
$\varpi_{S_2}$-configuration, $\bar{c}_2$, with the set of end
vertices $E_2\subseteq (B_1\cup \{x_1\})$
 and $S_2\subseteq A_2$. If $x_1\not\in E_2$, then let $\bar{c}_3$ be a blue $\varpi_{S_3}$-configuration with the set of end vertices
$E_3\subseteq\{x_1,y,z\}$ and $S_3\subseteq A_3$ where $y\in B_1$
and $z\in W\setminus (E_1\cup E_2)$. If $x_1\in E_2$, then let
$\bar{c}_3$ be a blue $\varpi_{S_3}$-configuration with the set of
end vertices $E_3\subseteq\{x_2,y,z\}$ and $S_3\subseteq A_3$
where $x_2\in E_2\setminus\{x_1\}$ and $\{y,z\}\subseteq
W\setminus(E_1\cup E_2)$. We continue this process to find the set
of $\{\bar{c}_1,\bar{c}_2,\ldots,\bar{c}_{q'}\}$ of
configurations. When this process terminate,
 %by combining these configurations
 we have the paths $\mathcal{P}_{l''}$ and
$\mathcal{P}_{l'}$ where $l''\geq l'\geq 0$ and $l''+l'=2q'$. Let
$x'',y''$ (resp. $x',y'$ if $l'> 0$) be the end vertices of
$\mathcal{P}_{l''}$ (resp. $\mathcal{P}_{l'}$) in $W$. Let
$T=V(\mathcal{K}^3_t)\setminus (V(\mathcal{P})\cup
V(\mathcal{P}_{l''})\cup V(\mathcal{P}_{l'}))$. Clearly
$|T|=\lfloor\frac{m+1}{2}\rfloor+1-(q'+i)$ where $i=1$ if $l'=0$
and $i=2$ if $l'>0$. Assume $m=4k+r$ for some $r$, $0\leq r\leq
3$. We have the following cases.

\bigskip
\noindent \textbf{Case 1. }$r=0$.

\medskip
Since $q\geq 2k-1$, we have $2q'\geq m-2$. On the other hand,
$|W|=\lfloor\frac{m+1}{2}\rfloor+1$ and so $2q'\leq m$. If $2q'=m$,
then $l'=0$ and so $\mathcal{P}_{l''=m}$ is a blue path. Now we
may assume that $2q'=m-2$, and one can easily check that the
vertices $\{v_{2n-3},v_{2n-2},v_{2n-1}\}$ are not used in
$\mathcal{P}_{l''}\cup\mathcal{P}_{l'}$. First let $l'=0$. Then
$|T|=1$ and we may assume $T=\{u\}$. Now using the maximality of
$\mathcal{P}$ and the
  fact that $\mathcal{C}^3_{n}\nsubseteq \mathcal{F}_{red}$, $\mathcal{P}_{l''}\cup\{\{v_{2n-2},y'',u\},
  \{v_{2n-1},u,v_1\}\}$ is a blue $\mathcal{P}^3_{m}$.
  For  $l'> 0$,
   $\mathcal{P}_{l''}\cup\{\{v_{2n-2},y'',x'\}\}\cup\mathcal{P}_{l'}\cup\{\{v_{2n-1},y',v_1\}\}$
   is a blue $\mathcal{P}^3_{m}$.

\bigskip
\noindent \textbf{Case 2. }$r=1$.

\medskip
Since $\mid W\mid=\lfloor \frac{m+1}{2}\rfloor+1$, $2q'\leq m+1$
and if the equality holds, then $l'=0$. On the other hand, $q\geq
2k$ and so $2q'\geq m-1$. Hence $2q'\in \{m+1,m-1\}$. If
$2q'=m+1$, then $l'=0$ and there is a blue $\mathcal{P}^3_{m+1}$.
Now let $2q'=m-1$. If $l'=0$, then $|T|=1$, so $T=\{u\}$ and hence
$\mathcal{P}_{l''}\cup\{\{v_1,u,y''\}\}$ is a blue
$\mathcal{P}^3_{m}$. If $l'>0$, then $\mathcal{P}_{l''}\cup
\{\{v_1,y'',x'\}\}\cup \mathcal{P}_{l'}$ is a blue
$\mathcal{P}^3_{m}$.

\bigskip
\noindent \textbf{Case 3. }$r=2$.

\medskip
Using an argument similar to the case 1, we have $2q'\in
\{m,m-2\}$ and if $2q'=m$, then $l'=0$ and we have a blue
$\mathcal{P}_{l''=m}$. Again by an argument similar to the case 1
we have a blue $\mathcal{P}^3_m$.

\bigskip
\noindent \textbf{Case 4. }$r=3$.

\medskip
In this case, partition  $V(\mathcal{P})\setminus \{v_1,v_2\}$
into $\lfloor\frac{2n-3}{5}\rfloor$ classes of size five and
possibly one class of size at most four. Then we repeat the
mentioned process in the first of the proof to find blue paths
$\mathcal{P}_{l''}$ and $\mathcal{P}_{l'}$ with $l''\geq
l'\geq 0$ and $l''+l'=2q'$. Again using a similar argument in case 1, we have $2q'\in
\{m+1,m-1,m-3\}$. If $2q'=m+1$, then we have $l'=0$ and so there is a
blue $\mathcal{P}^3_{m+1}$. For $2q'=m-1$, the assertion holds by
an argument similar to the case 2. Now let $2q'=m-3$. If $l'=0$,
then $|T|=2$, so $T=\{u,v\}$ and hence $\mathcal{P}_{l''}\cup
\{\{v_{2n-2},v_2,y''\}, \{v_{2n-2},v,u\}, \{u,v_1,v_{2n-1}\}\}$
is a blue $\mathcal{P}^3_{m}$( note that
$\{v_{2n-3},v_{2n-2},v_{2n-1}\}\cap V(\mathcal{P}_{l''})=\emptyset
$). If $l'>0$, then $|T|=1$, so $T=\{u\}$ and hence
$\mathcal{P}_{l''}\cup \{\{v_{2n-2},v_2,y''\},
\{v_{2n-2},x',u\}\}\cup \mathcal{P}_{l'}\cup \{\{y',v_1,v_{2n-1}\}\}$ is
a blue $\mathcal{P}^3_{m}$ and the proof is completed.
 $\hfill\blacksquare$

%%%%%%%%%%%%%%%%%%%%%%%%%%%%%%%%%%%%%%%%%%%%%%%%%%%%%%%%%%%%%%%%%%%%%%%%%%%%%%%%%%%%%%%%%%%%%%%%%%%%%%
%%%%%%%%%%%%%%%%%%%%%%%%%%%%%%%%%%%%%%%%%%%%%%%%%%%%%%%%%%%%%%%%%%%%%%%%%%%%%%%%%%%%%%%%%%%%%%%%%%%%%%

\bigskip
\begin{theorem}\label{main theorem}
For every $n\geq \Big\lfloor\frac{5m}{4}\Big\rfloor$, $$R(\mathcal{P}^3_n,\mathcal{P}^3_m)=2n+\Big\lfloor\frac{m+1}{2}\Big\rfloor.$$
\end{theorem}
%%%%%%%%%%%%%%%%%%%%%%%%%%%%
%%%%%%%%%%%%%%%%%%%%%%%%%%%%
\noindent\textbf{Proof. }We prove the theorem by induction on
$m+n$.  The proof of the case  $m=n=1$ is trivial. Suppose that for
$m^{\prime}+n^{\prime}< m+n$ with $n'\geq
\lfloor\frac{5m'}{4}\rfloor$, $R(\mathcal{P}^3_{n^{\prime}},
\mathcal{P}^3_{m^{\prime}})=2n^{\prime}+\Big\lfloor\frac{m^{\prime}+1}{2}\Big\rfloor.$
Now, let $n\geq \Big\lfloor\frac{5m}{4}\Big\rfloor$ and let
$\mathcal{K}^3_{2n+\lfloor\frac{m+1}{2}\rfloor}$ be 2-edge colored
red and blue. We may assume there is no red copy of
$\mathcal{P}^3_n$ and no blue copy of $\mathcal{P}^3_m$.  Consider
the following cases.

\bigskip
\noindent \textbf{Case 1. } $n= \Big\lfloor\frac{5m}{4}\Big\rfloor$.

\medskip
Since $R(\mathcal{P}^3_{n-1},\mathcal{P}^3_{m-1})=
2(n-1)+\Big\lfloor\frac{m}{2}\Big\rfloor<
2n+\Big\lfloor\frac{m+1}{2}\Big\rfloor$ by induction hypothesis,
then either there is a $\mathcal{P}^3_{n-1}\subseteq
\mathcal{F}_{red}$ or a $\mathcal{P}^3_{m-1}\subseteq
\mathcal{F}_{blue}$. If we have a red copy of
$\mathcal{P}^3_{n-1}$, then by Lemma \ref {a} we have a
$\mathcal{P}^3_{m}\subseteq \mathcal{F}_{blue}$. Now assume that
there is a blue copy of  $\mathcal{P}^3_{m-1}$.
%then Since $n-1= \Big\lfloor\frac{5(m-1)}{4}\Big\rfloor$,
 Lemma \ref{pm-1 implies pn-1} implies that $\mathcal{P}^3_{n-1}\subseteq
\mathcal{F}_{red}$ and using Lemma \ref {a} we have
$\mathcal{P}^3_{m}\subseteq \mathcal{F}_{blue}$, a contradiction.

\bigskip
\noindent \textbf{Case 2. }$n> \Big\lfloor\frac{5m}{4}\Big\rfloor$.

\medskip
 In this case, $n-1\geq \Big\lfloor\frac{5m}{4}\Big\rfloor$ and since $R(\mathcal{P}^3_{n-1},\mathcal{P}^3_{m})=
2(n-1)+\Big\lfloor\frac{m+1}{2}\Big\rfloor<
2n+\Big\lfloor\frac{m+1}{2}\Big\rfloor$, by induction hypothesis
we have a $\mathcal{P}^3_{n-1}\subseteq \mathcal{F}_{red}$.
 Using  Lemma \ref {a} we have a
$\mathcal{P}^3_{m}\subseteq \mathcal{F}_{blue}$ and it completes
the proof.
 $\hfill\blacksquare$

\section{\normalsize Acknowledgments}

\emph{\emph{ The authors appreciate the discussions on
the subject of this paper with A. Gy\'arf\'as.  } }

%%%%%%%%%%%%%%%%%%%%%%%%%%%%%%%%%%%%%%%%%%%%%%%%%%%%%%%%%%%%%%%%%%%%%%%%%%%%%%%%%%%%%%%%%%%%%%%%%%%%%%%%%%%
%%%%%%%%%%%%%%%%%%%%%%%%%%%%%%%%%%%%%%%%%%%%%%%%%%%%%%%%%%%%%%%%%%%%%%%%%%%%%%%%%%%%%%%%%%%%%%%%%%%%%%%%%%
%%%%%%%%%%%%%%%%%%%%%%%%%%%%%%%%%%%%%%%%%%%%%%%%%%%%%%%%%%%%%%%%%%%%%%%%%%%%%%%%%%%%%%%%%%%%%%%%%%%%%%%%%%%
%%%%%%%%%%%%%%%%%%%%%%%%%%%%%%%%%%%%%%%%%%%%%%%%%%%%%%%%%%%%%%%%%%%%%%%%%%%%%%%%%%%%%%%%%%%%%%%%%%%%%%%%%%%

\end{document}